\documentclass{amsart}

\newtheorem{theorem}{Theorem}[section]

\theoremstyle{definition}

\theoremstyle{remark}

\numberwithin{equation}{section}

    \newcommand{\vsc}{\vspace{6pt}}
  \newcommand{\R}{{I\!\! R}}
  \newcommand{\CO}{\mathcal{O}}
   \newcommand{\tv}{\hbox{Tot.Var.}}
   \newcommand{\N}{\mathcal{N}}
   \newcommand{\BV}{{\bf BV \/}}
   \newcommand{\C}{\mathcal{C}}
   \newcommand{\ep}{\epsilon}
   \newcommand{\CR}{\mathcal{R}}
   \newcommand{\D}{\mathcal{D}}
   \newcommand{\Z}{Z\!\!\!Z}
   \newcommand{\be}{\begin{equation}}
   \newcommand{\ee}{\end{equation}}
    \newcommand{\rrn}{{I\!\! R^n}}
    \newcommand{\lu}{{{\bf L}^1}}
    \newcommand{\li}{{{\bf L}^\infty}}
    \newcommand{\G}{\mathcal{G}}
    \newcommand{\J}{\mathcal{J}}
    \newcommand{\CS}{\mathcal{S}}
    \newcommand{\CZ}{\mathcal{Z}}
    \newcommand{\CP}{\mathcal{P}}
    \newcommand{\A}{\mathcal{A}}
    \newcommand{\T}{\mathcal{T}}

\begin{document}

\title[Decay of Positive Waves for $n \times n$ Balance Laws]{Decay of Positive 
Waves for $n \times n$ Hyperbolic Systems of Balance Laws}

\author{Paola Goatin}
\address{ Centre de Math\'ematiques Appliqu\'ees,
Ecole Polytechnique, 91128 Palaiseau Cedex, France.}
\email{goatin@cmap.polytechnique.fr}
\thanks{The authors were partially supported respectively by the EC-Marie Curie
Individual Fellowship \#HPMF-CT-2000-00930 and EEC grants
\#ERBFMRXCT970157 \& \#HPRN-CT-2002-00282.}

\author{Laurent Gosse}
\address{Istituto per le Applicazioni del Calcolo (sezione di Bari), via G. Amendola, 122/I, 
70126 BARI, Italy } \email{l.gosse@area.ba.cnr.it}


\date{February 2003}



\begin{abstract}
We prove Ole\u \i nik-type decay estimates for entropy solutions of $n\times n$ strictly hyperbolic
systems of balance laws built out of a wave-front tracking procedure inside which the source
term is treated as a nonconservative product localized on a discrete lattice.

\end{abstract}

\maketitle

\section{Introduction}

A classical result proved by Ole\u \i nik \cite{ol} for strictly convex scalar 
conservation laws in one space dimension shows that the density of positive waves 
decays in time like $\CO (1/t)$, see also \cite{hoff}. More precisely, if we
consider the scalar equation \begin{equation} 
u_t + f(u)_x = 0, \qquad u(t=0,.)=u_o \in \li(\R),
\end{equation}
with $f''(u)\geq \kappa >0$, we have that every entropy-admissible 
solution satisfies
\begin{equation} 
u(t,y)-u(t,x)\leq \frac{y-x}{\kappa t} \qquad\hbox{for all $t>0$, $x<y$}.
\label{oleinik} \end{equation}
and therefore has locally bounded variation (see also \cite{daf}, Theorem 11.2.2).
Conversely, if $u=u(t,x)$ is a weak solution satisfying $(\ref{oleinik})$, then $u$ is entropy 
admissible. 

The same estimate as in $(\ref{oleinik})$ has been recovered for the Riemann coordinates 
of a particular $2\times 2$ system \cite{brco} and for $n\times n$ genuinely nonlinear systems 
belonging to the Temple class \cite{brgo2}.
However, one cannot expect such a result to remain valid for general $n\times n$ systems, even assuming
all characteristic fields being genuinely nonlinear. Indeed, interactions among existing shocks
may generate rarefactions as time increases. Decay estimates must therefore take into account
for the generation of new positive waves due to interactions. Results in this direction were
proved by Liu \cite{liu} in the case of approximate solutions constructed by
Glimm's scheme, \cite{glimm}, and by Bressan and Colombo \cite{brcol} for exact
solutions obtained as limits of front tracking appoximations for $n \times n$
homogeneous systems, \cite{nils}. This in turn yields uniqueness of solutions
satisfying the Ole\u \i nik entropy condition \cite{brgo1, brgo2,paola}. From the
point of view of practical applications, such one-sided estimates
are very useful for instance in the context of multiphase Geometric Optics
computations, \cite{JCP}, or local error estimates, \cite{tad} and asymptotic
behaviour of entropy solutions.

In this paper we are interested in extending Ole\u \i nik-type estimates on
positive waves  to quasilinear systems of balance laws. More precisely, we
shall deal with the Cauchy problem for the following $n\times n$ system of
equations  \begin{equation} 
u_t + f(u)_x = g(x,u), \qquad x\in\R, t>0, \label{bl}
\end{equation}
endowed with a (suitably small) initial data $u_o \in \lu\cap\BV (\R;\rrn)$.
Here $u(t,x)\in\rrn$ is the unknown function, $f:\Omega\to\rrn$ is a smooth
${\mathcal C}^2$ vector field defined on an open neighborhood $\Omega$ of the
origin in $\rrn$. We will assume that the system $(\ref{bl})$ is strictly
hyperbolic, with each  characteristic field either genuinely nonlinear or
linearly degenerate in the  sense of Lax \cite{Lax}. Moreover we assume the
following  Caratheodory-type conditions for the source term $g$:
\begin{itemize}
\item[(P1)] $g:\R\times\Omega\to\rrn$ is measurable w.r.t. $x$, for any $u\in\Omega$,
            and is $\C^2$ w.r.t. $u$, for any $x\in\R$;
\item[(P2)] $\| g(x,\cdot)\|_{\C^2}$ is bounded over $\Omega$, uniformly in $x\in\R$;
\item[(P3)] there exists a function $\omega\in\lu(\R)\cap\li(\R)$ such that
            $|g(x,u)|\leq\omega(x)$, $\|\nabla_u g(x,u)\|\leq\omega(x)$
            for all $(x,u)\in\R\times\Omega$.
\end{itemize}
In addition, we require that a {\it non-resonance} condition holds, that is, the characteristic
speeds are bounded away from zero: for some $p\in \{1,\ldots,n\}$ and some $c>0$ one has
\begin{equation}
\begin{cases}
\lambda_i (u)\leq -c & \mbox{ if $i\leq p$}, \\
\lambda_i (u)\geq c  & \mbox{ if $i> p$},
\end{cases} 
\label{nonres} \end{equation}
for all $u\in\Omega$, where $\lambda_i (u)$ denote the eigenvalues of the Jacobian matrix $Df(u)$.

Under these assumptions, it was proved in \cite{agg} that there exists a family of entropy 
weak solutions to $(\ref{bl})$ continuously depending on the initial data. More precisely,
if the $\lu$-norm of $\omega$ is small enough, there exist a closed domain
$\D\subset\lu (\R;\rrn)$ of functions with sufficiently small total variation, a constant $L$
and a unique semigroup $P:[0,+\infty)\times\D\to\D$ with the properties:
\begin{itemize}
\item[(i)] For all $u, v\in\D$ and $t,s\geq 0$ one has 
           $\|P_s u -P_t v\|_{\lu} \leq L \big( |t-s|+\|u-v\|_{\lu}\big)$.
\item[(ii)] For all $u_o \in\D$ the function $u(t,\cdot)=P_t u_o$ is a weak entropy
            solution of the Cauchy problem $(\ref{bl})$, $u(t=0,.)=u_o$.
\end{itemize}
Under the above assumptions we aim at showing that, for genuinely nonlinear
characteristic fields,  an Ole\u \i nik type estimate on the decay of positive
waves holds, which takes into account  not only new waves generated by
interactions but even the contribution of the source term. A careful statement
of these results requires some notations, \cite{bressan,brcol,brgo1}.

As usual, let $A(u)=Df(u)$ be the Jacobian matrix of $f$, and call $\lambda_i(u)$,
$l_i(u)$, $r_i(u)$ respectively the eigenvalues and the left and right eigenvctors of $A(u)$.
Let $u:\R\to\Omega$ have bounded variation in $x$ and satisfy (\ref{bl}) with $g \equiv 0$. 
The distributional
derivative  $\mu\doteq D_x u$ is a vector measure. For $i=1,\ldots,n$ we can
now define $\mu^i$ as $$     
\int \phi d\mu^i =\int \phi \tilde l_i \cdot D_x u, \qquad \phi\in\C^0_c,
$$
where $\tilde l_i (x) = l_i(u(x))$ at points where $u$ is continuous, while
$\tilde l_i (x_\alpha)$ is some vector which satisfies
$$
|\tilde l_i (x_\alpha) -l_i (u(x_\alpha))| = \CO(1)\cdot |u(x_\alpha+)-u(x_\alpha-)|
$$
$$
\tilde l_i (x_\alpha)\cdot\big(u(x_\alpha+)-u(x_\alpha-)\big)=\sigma_\alpha^i,
$$
where with $\sigma_\alpha^i$ we denote the strength of the $i$-th wave generated by 
the resolution of the corresponding discontinuity in $x_\alpha$.

We denote by $\mu^{i+}$, $\mu^{i-}$ the positive and negative parts of $\mu$, then we have
$$
\mu^i = \mu^{i+}-\mu^{i-}, \qquad |\mu^i| = \mu^{i+}+\mu^{i-}.
$$
The {\it total strength of waves} in $u$ is defined as
\begin{equation} \label{forza}
{\bf V}(u)\doteq \sum_{i=1}^n {\bf V}_i (u), \qquad {\bf V}_i (u)\doteq |\mu^i| (\R),
\end{equation}
while the {\it interaction potential} is defined in terms of product measures on $\R^2$:
\begin{equation} \label{potenziale}
{\bf Q}(u) \doteq \sum_{i<j}\big( |\mu^j|\times |\mu^i|\big) \big( \{ (x,y): x<y \}\big)
           + \sum_{i\in\G\N}\big( \mu^{i-}\times |\mu^i|\big) \big( \{ (x,y): x\not= y \}\big),
\end{equation}
where $\G\N$ denotes the set of genuinely nonlinear families.

Now we are ready to state our main result in the case $g \not \equiv 0$:
\begin{theorem} \label{teorema}
{\rm\bf (Decay of positive waves)} Let the system $(\ref{bl})$ 
be strictly hyperbolic and let the $i$-th characteristic field be genuinely 
non-linear. Then there exists a constant $C$ depending solely on $f$
such that for every $0\leq s<t$ and every solution $u$ with small total
variation obtained as limit of wave-front tracking approximations, the measure
$\mu^{i+}_t$ of $i$-waves  in $u(t,\cdot)$ satisfies
\begin{equation} \label{olly}
\mu^{i+}_t (J)\leq C\cdot \left(\frac{meas(J)}{t-s} + {\bf Q}(s)-{\bf
Q}(t)  +{\bf V}(u_o)\cdot\|\omega\|_\lu \right)
\end{equation}
for every Borel set $J\subset\R$.
\end{theorem}

Of course, we tacitly assume (P1)--(P3) and (\ref{nonres}) throughout all the
text.

\section{Wave-front tracking with zero-waves}

In this section we briefly recall the construction of wave-front tracking
approximations as stated in \cite{agg}.
We start with the definition of the $h$-Riemann solver. For small $h>0$
we introduce the map
$$
\Phi_h(x_o,u)\doteq f^{-1}\left[ f(u)+\int_0^h g(x_o+s,u)ds \right]
$$
(note that $f$ is invertible due to $(\ref{nonres})$), that approximates 
the flow of the stationary equation associated to $(\ref{bl})$.
Consider now the Riemann problem with initial states 
\begin{equation}
u(0,x)=\left\{ 
\begin{array}{ll}
u_l  &\textrm{if $x<x_o$}, \\
u_r &\textrm{if $x>x_o$}.
\end{array} \right. \label{riemann}
\end{equation}
To locally render the source term's effects, a stationary discontinuity 
is introduced along the line $x=x_o$, that is, a wave whose speed is 
equal to zero; it will be referred to as a {\bf zero-wave}.
An $h$-Riemann solver for $(\ref{bl})-(\ref{riemann})$ has been defined in
\cite{agg} as a self-similar function $u(t,x)=R_h((x-x_o)/t;u_l,u_r)$ as
follows: %
\begin{itemize}
\item[(a)] there exist two states $u^-$, $u^+$ which satisfy 
           $u^+ = \Phi_h (x_o,u^-)$;
\item[(b)] $u(t,x)$ coincides, on the set $\{t\geq 0,x<x_o\}$, with the solution to 
           the homogeneous Riemann problem with initial values $u_l$, $u^-$
	   and, on the set $\{t\geq 0,x>x_o\}$, with the solution to the
homogeneous            Riemann problem with initial values $u^+$, $u_r$;
\item[(c)] the Riemann problem between $u_l$ and $u^-$ is solved only by waves with 
           negative speed (i.e. of the families $1,\ldots,p$);
\item[(d)] the Riemann problem between $u^+$ and $u_r$ is solved only by waves with 
           positive speed (i.e. of the families $p+1,\ldots,n$).
\end{itemize}
This clearly shares a lot of common features with the nonconservative Riemann
problems studied in \cite{lft}. Let now $\ep,h>0$ be given: an
$\ep,h$-approximate solution of $(\ref{bl})$ is constructed as follows.
First of all, the source term is localized by means of a Dirac comb along zero-waves 
located on the lattice $x=jh$, $j\in(-{1\over h\ep},{1\over h\ep})\cap\Z$:
\begin{equation}
u_t +f(u)_x = h \sum_j g(x,u).\delta (x-jh), \label{dirac}
\end{equation}
where $\delta$ stands for the Dirac measure concentrated on $x=0$.

Given the initial data $u_o$, we deduce a piecewise constant approximation
$u(0,\cdot)$ without increasing its {\bf BV}-norm and $u(t,x)$ is
constructed, for small $t$, by applying the $h$-Riemann solver at every point
$x=jh$, and by solving the remaining discontinuities in  $u(0,\cdot)$ using a
classical homogeneous Riemann solver (rarefaction waves are discretized 
following \cite{bressan}: for a fixed small parameter $\nu$, each rarefaction
of size $\sigma$ is divided, at its starting time, into
$m=\left[\sigma\over\nu\right]+1$ wave fronts of size $\sigma/m\leq\nu$).

At every interaction point, a new Riemann problem arises. Notice that because of their
null speed, zero-waves cannot interact among each other. In order to keep
finite the  total number of wave-fronts, two distinct procedures are used for
solving a Riemann problem: an accurate method, which possibly creates
several new fronts, and a simplified method , which minimizes the number of
new wave-fronts. For a detailed description, as well as the proof of the
stability of the algorithm, we refer the reader to  \cite{agg, bressan}.

The approximate solution can have four types of jumps: shocks (or contact discontinuities),
rarefaction fronts, non-physical waves and zero-waves: $\J
=\CS\cup\CR\cup\N\CP\cup\CZ$. A priori bounds on the functions $u^{\ep,h}$
are obtained by modifying slightly the Glimm's functionals 
\cite{glimm} in order to keep track of the zero-waves,
\begin{equation}\label{tv}
V(t) = \sum_{\alpha\in\J} |\sigma_\alpha| = \sum_{\alpha\in\CS\cup\CR\cup\N\CP\cup\CZ} 
|\sigma_\alpha|, \end{equation}
\begin{equation}\label{interaction}
Q(u(t)) = \sum_{\alpha,\beta\in\tilde\A} |\sigma_\alpha \sigma_\beta| \leq V(t)^2,
\end{equation}
measuring respectively the total wave strengths and the interaction potential in $u(t,\cdot)$.
In particular, if $\alpha\in\CZ$ then the strength of the wave located in 
$x_\alpha = j_\alpha h$ can be measured by means of  
\begin{equation}
\sigma_\alpha =\int_0^h \omega(j_\alpha h+s)ds. \label{forzazero}
\end{equation}
In $(\ref{interaction})$ we have denoted by $\tilde\A$ an extended set of {\it
approaching waves}.  As usual we call $k_\alpha$ the family of the front
located at $x_\alpha$, with size $\sigma_\alpha$. More precisely, a couple of 
wave-fronts of families $k_\alpha$, $k_\beta$, 
located at $x_\alpha < x_\beta$,  belongs to $\tilde\A$ in any of the following cases:
\begin{itemize}
\item[--]  if none of the two is a zero-wave, either $k_\alpha >k_\beta$, or
else $k_\alpha = k_\beta$    and at least one of them is a genuinely nonlinear
shock, \cite{glimm}; %
\item[--]  if $\alpha$ is a zero-wave and $\beta$ is a physical one,
$k_\beta\leq p$; %
\item[--]  if $\beta$ is a zero-wave and $\alpha$ is a physical one,
$k_\alpha > p$. \end{itemize}
Notice that for some $C>1$ there holds
$$
{1\over C}\Big[ \|\omega\|_{\lu(I)}+\tv u(t,\cdot)\Big] \leq V(t)
\leq C \Big[ \|\omega\|_{\lu(I)}+\tv u(t,\cdot)\Big],
$$
where the interval $I$ is defined by
$$
I\doteq \bigcup_{\alpha\in\CZ} [j_\alpha h, (j_\alpha +1)h]
 =\bigcup_{j\in\left(-{1\over \ep h},{1\over \ep h}\right)\cap\Z} [j h, (j+1)h].
$$
Passing to the limit as $\ep\to 0$, $h>0$ fixed, one has:
\begin{itemize}
\item[(i)] the total variation of $u^{\ep,h} (t,\cdot)$ remains uniformly bounded;
\item[(ii)] the maximum size of rarefaction fronts approaches zero;
\item[(iii)] the total strength of all non-physical waves approaches zero.
\end{itemize}
By $(i)$, Helly's theorem guarantees the existence of a subsequence strongly convergent
in $\lu_{loc}$. By $(ii)$ and $(iii)$, this limit provides a weak solution to $(\ref{dirac})$
in agreement with nonconservative theories, \cite{lft}.
At this stage, we can extract again a
subsequence $u^{h_i}$ which converges to some function $u$ in $\lu_{loc}$
and solves $(\ref{bl})$ in the usual weak distributional sense.

\vsc

In the sequel, we will need a semicontinuity property of Glimm's functionals.
For $h>0$ fixed, let the total strength of waves ${\bf V}_h$ and the
interaction potential ${\bf Q}_h$ be as in \cite{agg}, Section 4.1, that is,
they are defined as in $(\ref{forza})$-$(\ref{potenziale})$, but including the
zero-waves. We have \begin{equation} \label{hforza}
{\bf V}_h (u) = {\bf V} (u) + \| \omega\|_\lu,
\end{equation}
\begin{equation} \label{hpotenziale}
{\bf Q} (u)\leq {\bf Q}_h \leq {\bf Q} (u) +\| \omega\|_\lu \cdot {\bf V} (u).
\end{equation}
Proceeding as in \cite{bressan} we recover
the lower semicontinuity of the functionals ${\bf Q}_h$ and 
${\bf \Upsilon}_h (u)\dot={\bf V}_h (u)+C_o {\bf Q}_h (u)$, $C_o>0$, on a domain $\D$
of the form $$
\D\dot= \Big\{ u\in\lu \cap {\bf BV}(\R;\rrn), {\bf \Upsilon}_h
(u)\leq\gamma\Big\}, \qquad \gamma \mbox{ small enough,} $$
(see \cite{bressan}, Theorem 10.1):
\begin{theorem} \label{lsc}
{\rm\bf (Lower semicontinuity of the Glimm functionals)} There exists a choice
of the constants $C_o, \gamma >0$ such that, if ${\bf \Upsilon}_h (u)={\bf V}_h
(u)+C_o {\bf Q}_h (u) <\gamma$, then for any sequence of
functions $u_\nu \in\D$, $u_\nu \to u$ in $\lu$ as $\nu\to\infty$, one has
\begin{eqnarray}
{\bf Q}_h (u)\leq \liminf_{\nu\to\infty} {\bf Q}_h (u_\nu), \label{qlsc}\\
{\bf \Upsilon}_h (u) \leq\liminf_{\nu\to\infty} {\bf \Upsilon}_h (u_\nu). \label{glsc}
\end{eqnarray}
Moreover, for every finite union of open intervals $J=I_1 \cup\ldots\cup I_S$
there holds: \begin{equation} \label{mulsc}
\mu^{i\pm}(J)+C_o {\bf Q}_h (u)\leq \liminf_{\nu\to\infty}\big( \mu^{i\pm}_\nu(J)+C_o {\bf Q}_h (u_\nu)\big),
\qquad i=1,\ldots,n.
\end{equation}
\end{theorem}

\section{Proof of Theorem \ref{teorema}}

{\bf I. } By Lipschitz continuous dependence of the trajectories it is not
restrictive to assume $s=0, t=T$. We will consider a particular converging
sequence $u^{\nu,h}$ of $\ep_\nu,h$-approximate solutions with the following
properties: %
\begin{itemize}
\item[(i)] each rarefaction front $x_\alpha$ travels with the characteristic speed of the state 
           on the right:
           $$
           \dot x_\alpha =  \lambda_{k_\alpha} (u(x_\alpha+));
           $$
\item[(ii)] each shock $x_\alpha$ travels with a speed strictly contained between 
            the right and the left characteristic speeds:
            $$
            \lambda_{k_\alpha} (u(x_\alpha+))<\dot x_\alpha <\lambda_{k_\alpha} (u(x_\alpha-));
            $$
\item[(iii)] calling $N_\nu$ the number of jumps in $u_o^{\nu,h}=u^{\nu,h}(0,\cdot)$, as $\nu\to\infty$
             one has
             \begin{equation} \label{enlim}
             \ep_\nu \to 0 \qquad  \ep_\nu N_\nu \to 0 ;
             \end{equation}
\item[(iv)] the interaction potential satisfies
            \begin{equation} \label{inpot}
            {\bf Q}_h(u^{\nu,h}(0,\cdot)) \to {\bf Q}_h (u_o) \qquad  \hbox{as $\nu\to\infty$}.
            \end{equation}
\end{itemize}
Such a sequence can be constructed as explained in \cite{bressan}, proof of Lemma 10.2 (p.205).

Let $u=u(t,x)$ be a piecewise constant $\ep,h$-approximate solution
constructed via front-tracking approximation (we shall drop the $\ep,h$
superscripts since there is no ambiguity). As usual, by {\it (generalized) 
$i$-characteristic} we mean an absolutely continuous curve $x=x(t)$ such that,
\cite{daf}, $$ \dot x (t)\in [\lambda_i (u(t,x+)), \lambda_i (u(t,x-))]
\qquad\hbox{a.e. $t\geq 0$} $$
By $t\mapsto y^i(t;\bar x)$ we denote the minimal $i$-characteristic passing through $\bar x$ at time $T$.
Because of $(\ref{nonres})$ the presence of zero-waves does not affect the usual construction.

Let now $I\doteq [a,b[$ be any half-open interval, and define
$$
I(t)\doteq [ y^i(t;a), y^i(t;b)[\doteq [a(t),b(t)[.
$$ 
We seek an estimate of the amount of positive $i$-waves in the the approximate solution $u(T,\cdot)$
contained in $I$. We recall that $k_\alpha$ stands for the family of the front
located at $x_\alpha$, with size $\sigma_\alpha$. For a genuinely nonlinear family, 
the size of the jump can be measured like
\begin{equation}
\sigma_\alpha\doteq \lambda_{k_\alpha} (u(x_\alpha+))-\lambda_{k_\alpha}
(u(x_\alpha-)),  \label{forzagnl} \end{equation}
while the size of a zero-wave is still given by $(\ref{forzazero})$. Define
$$
m(t)\doteq b(t)-a(t).
$$
By $(\ref{forzagnl})$ and the Lipschitz continuity of the map $u\mapsto\lambda_i(u)$ we deduce that
\begin{eqnarray}
\dot m(t) &=& \lambda_i\big(u(t,b(t))\big)-\lambda_i\big(u(t,a(t))\big) 
\nonumber \\           &=& M(t) +\CO (1) (\ep+K(t)) \label{eqbase}  
\end{eqnarray}
for a.e. $t$. The Landau symbol stands for a quantity whose modulus is
uniformly bounded. We see that  $$
M(t)\doteq \sum_{k_\alpha =i,x_\alpha\in I(t)} \sigma_\alpha
=\mu_t^i\Big(I(t)\Big)
$$
is the total amount of (signed) $i$-waves in $u(t,\cdot)$ contained in $I(t)$, while
$$
K(t)\doteq \sum_{k_\alpha \not= i,x_\alpha\in I(t)} |\sigma_\alpha|
=\sum_{k \not= i}\big|\mu_t^k\big|\Big(I(t)\Big)+\int_{I(t)}\omega(x).dx
$$
stands for the total strength of waves of families $\not= i$ inside
$I(t)$,  zero-waves included. To estimate the contribution of the term $K(t)$
in $(\ref{eqbase})$ we introduce 
$$
\Phi (t)\doteq\sum_{k_\alpha\not= i}\phi_{k_\alpha}(t,x_\alpha(t))\cdot |\sigma_\alpha|
\leq {\bf V}_h(u(t)),
$$
where
$$
\phi_j (t,x)\doteq \left\{
  \begin{array}{ll}
          1 &\textrm{if  $x<a(t)$} \\
          {{b(t)-x}\over{m(t)}} &\textrm{if $x\in [a(t),b(t)[$}\\
          0 &\textrm{if $x\geq b(t)$}
  \end{array}\right.\,
$$
or
$$
\phi_j (t,x)\doteq \left\{
  \begin{array}{ll}   
          0 &\textrm{if $x<a(t)$}  \\
          {{x-a(t)}\over{m(t)}} &\textrm{if $x\in [a(t),b(t)[$}  \\
          1 &\textrm{if $x\geq b(t)$}
  \end{array}\right.
$$          
in the cases $j<i$ or $j>i$ respectively. Roughly speaking, $\Phi(t)$ represents the cumulated 
strength of the waves which do {\bf not} approach the interval $I(t)$.
By strict hyperbolicity, we can expect it to grow with time. Observe that
$\Phi$ is piecewise Lipschitz continuous with  a finite number of
discontinuities occurring at interaction times, where it may decrease at most
of: \begin{equation}
\Phi (\tau+)-\Phi (\tau-)=\CO (1) [Q(u(\tau-))-Q(u(\tau+))]. \label{salti}
\end{equation}
Following \cite{bressan}, we assume that there holds for
some $c_o>0$, \begin{equation}\label{nana} |\lambda_i(u)-\lambda_i(v)|\leq
c_o,\qquad |\lambda_i(u)-\lambda_{k_\alpha}(v)|\geq 2c_o, \end{equation}
for every couple of states $u,v$ and every $k_\alpha\not= i$. Outside
interaction times $\Phi$ is non-decreasing; indeed, we have 
\begin{eqnarray*}
\dot \Phi (t)&=&\sum_{k_\alpha\not= i} |\sigma_\alpha|\cdot {d\over dt}\phi_{k_\alpha}(t,x_\alpha(t))\\
             &=&\sum_{k_\alpha<i,x_\alpha\in I(t)} |\sigma_\alpha|\cdot
                \left({\dot b -\dot x_\alpha \over m}-{(b-x_\alpha)\dot m \over m^2}\right) \\
             & &+\sum_{k_\alpha>i,x_\alpha\in I(t)} |\sigma_\alpha|\cdot
                \left({\dot x_\alpha -\dot a \over m}-{(x_\alpha -a)\dot m \over m^2}\right) \\
             &\geq&\sum_{k_\alpha\not= i} |\sigma_\alpha|\cdot {c_o\over m(t)} \\ \end{eqnarray*} 
thanks to the system's strict hyperbolicity and
the non-resonance condition $(\ref{nonres})$. In particular, observe that
$c_o\leq c$ in (\ref{nana}). The above estimate yields the bound valid for all
but finitely many times $t$: \begin{equation}
K(t)\leq {1\over c_o} \dot\Phi (t) m(t). \label{kbound}
\end{equation}
We notice again, as in \cite{agg,agg'}, that there's a need for a completely
different theory in order to tackle resonant cases, see also \cite{ha}.

Concerning the term $M(t)$, observe that it can change only when an interaction
occurs within the interval $[a(t),b(t)]$. In this case, one has
$$
M(\tau+)-M(\tau-) =\CO (1) [Q(u(\tau-))-Q(u(\tau+))].
$$
This yields an estimate of the form
\begin{equation}
M(T)-M(t) =\CO (1) \sum_{\tau\in\T} [Q(u(\tau-))-Q(u(\tau+))] \label{mbound}
\end{equation}
where the summation extends over all times $\tau\in ]0,T]$ where an interaction occurs
inside $[a(\tau),b(\tau)]$. 
Inserting the estimates $(\ref{kbound})-(\ref{mbound})$ in $(\ref{eqbase})$ we obtain
\begin{equation}\label{ineq}
\dot m(t)+C\dot\Phi(t)m(t)\geq M(T)-C\left(\ep+\sum_{\tau\in\T}|\Delta Q(\tau)|\right),
\end{equation}
for some constant $C$ and a.e. $t$. We now observe that $m$ is a continuous,
piecewise linear function of $t$, and $\Phi$ is uniformly bounded. It can
decrease only at interaction times, where $(\ref{salti})$ holds. Hence its
total variation in uniformly bounded and for some  constant $K_o$ we have the
estimate %
\begin{equation}
\int_0^T \dot\Phi(t)dt\leq K_o. \label{intbound}
\end{equation}
As in \cite{bressan}, from $(\ref{ineq})$ we deduce the decay estimate:
\begin{equation} \label{basicin}
M(T)\leq 2e^{CK_o}\left({b-a\over
T}\right)+2C\ep+2C \sum_{\tau\in\T}|\Delta Q(\tau)|. \end{equation}
{\bf II. } Repeating the above process for any finite number $S$ of disjoint
half-open intervals $I_s\dot= [a_s,b_s[$ we obtain
\begin{equation} \label{sumin}
\sum_{s=1}^S M_s(T)\leq C' \left(\sum_{s=1}^S {b_s-a_s\over T}+S\ep+[Q(u(0))-Q(u(T))]\right),
\end{equation}
for some constant $C'$ independent of $S$ and of the particular $\ep,h$-approximate solution.
Here we have used the notation
$$
M_s (T)\doteq \sum_{k_\alpha=i, x_\alpha\in [a_s,b_s[} \sigma_\alpha =\mu_T^i([a_s,b_s[)
$$
to denote the sum of (signed) strength of all $i$-waves in $u(T,\cdot)$ contained in
the interval $[a_s,b_s[$. 

{\bf III. } Let us consider now any open interval $]a,b[$. Let $N$ be the
number of $i$-shocks of the first generation in the front-tracking
$\ep,h$-approximate solution $u$, as defined in \cite{bressan}, Chapter 7. We
can thus construct half-open intervals $I_s\dot= [a_s,b_s[$, 
$s=1,\ldots,S\leq N+1$, such that the following holds (see \cite{bressan},
Chapter 10): 
\begin{itemize}
\item Every $i$-rarefaction front in $u(T,\cdot)$ contained in $]a,b[$ falls inside one 
      of the intervals $I_s$.
\item No $i$-shock front of the first generation fall inside any of the intervals $I_s$.
\end{itemize}
Calling $\mu_T^{i+}$ the measure of positive $i$-waves in $u(T,\cdot)$, the above properties
imply
\begin{equation} \label{posmis}
\mu_T^{i+}(]a,b[)=\sum_s M_s(T)+\CO (1)\cdot [{\bf Q}_h (u(0))-{\bf Q}_h (u(T))]+\CO (1)\cdot r(\ep)
\end{equation}
where ${\bf Q}_h$ is the interaction potential introduced in \cite{agg}, Section 4.1.   
 
Indeed the only negative $i$-waves contained in $\bigcup_s I_s$ must have
generation order $\geq 2$, originating from interactions during the time
interval $]0,T]$. The total strength of these negative $i$-waves is bounded by
the decrease in the interaction potential $Q$. The last term on the right
hand side of $(\ref{posmis})$ tends to zero as $\ep$ does and comes from
the difference between $Q(u)$ and ${\bf Q}_h (u)$: in the latter there are no
non-physical fronts and all the countable $h$-Riemann problems are solved
(whereas before only zero-waves inside the interval $(-1/\ep,1/\ep)$ were
considered). 

Altogether, $(\ref{sumin})$ and $(\ref{posmis})$ yield
\begin{equation} 
\mu_T^{i+}(]a,b[) 
\leq C''\left({b-a \over T}+(N+1) r(\ep)+
                          [{\bf Q}_h (u(0))-{\bf Q}_h (u(T))]\right) \label{finale} 
\end{equation}
for some constant $C''$ independent of $\ep$.

{\bf IV. } For $h>0$ fixed, we now consider a sequence of
$\ep_\nu,h$-approximate solutions  satisfying the properties $(i)-(iv)$ stated
at the beginning of this section. It is clearly not restrictive to take
$C''\geq C_o$ in $(\ref{finale})$, where $C_o$ is the (big) constant in Theorem
\ref{lsc}. Using  $(\ref{qlsc})$, $(\ref{mulsc})$, $(\ref{finale})$,
$(\ref{enlim})$ and $(\ref{inpot})$ we obtain 
\begin{eqnarray*} 
& &\mu_T^{i+} (]a,b[) \\ 
& &\leq \liminf_{\nu\to\infty}
           \big( \mu_{\nu,T}^{i+} (]a,b[) +C''{\bf Q}_h(u_\nu(T))\big) -C''{\bf Q}_h(u(T)) \\
& &\leq C''\cdot\liminf_{\nu\to\infty}\left({b-a \over T}+(N_\nu +1) r(\ep)+{\bf Q}_h (u_\nu(0))
     \right) -C''{\bf Q}_h (u(T)) \\ 
& &\leq C''\cdot {b-a \over T}  +C''\Big[{\bf Q}_h (u(0))-{\bf Q}_h (u(T))\Big] \\
& &\leq C''\cdot {b-a \over T} +C''\Big[{\bf Q} (u(0))-{\bf
Q}(u(T)) +\|\omega\|_\lu \cdot {\bf V}(u(0))\Big]. \end{eqnarray*}
Since the last term is independent of $h$
this proves $(\ref{olly})$ in case $J$ is an open interval. The same arguments
can be used in the case where $J$ is a finite collection of open intervals.
Since $\mu^i$ is a bounded Radon measure, the estimate $(\ref{olly})$ holds for
every Borel set $J$ and we are done. $\Box$

\section{A short comment}

All these computations heavily rely on the restrictive assumption that the source term is
dominated by a function $\omega\in\lu(\R)$ (the $\li(\R)$ bound
being just a consequence of the smoothness of $g$).
In the case of a convex scalar law, \cite{CRAS},
$$
u_t+f(u)_x=g(u), \qquad u_o \in \lu \cap {\bf BV}(\R), 
$$
the stability of the approximation procedure for (\ref{bl}) can be obtained requiring only $\omega
\in \li(\R)$.
Hence one can follow the same canvas and (\ref{eqbase})
simplifies a lot thanks to the obvious bound $K(t) \leq \|\omega\|_{\li}
m(t)$. Thus (\ref{ineq}) boils down to 
$$
\dot m(t)+ \|\omega\|_{\li} m(t) \geq M(T),
$$
which leads in sharp contrast to an exponential bound of the type
($\mu^1=u_x$): \begin{equation}\label{milo_marat}
\mu^{1+}_t(J)\leq C e^{\|\omega\|_{\li} t}\cdot\frac{meas(J)}{t}, \qquad t>0.
\end{equation}
This agrees of course with simple computations, since in this context
Ole\u \i nik's estimates can be derived from the Riccati
differential equation: $$
z_t-g'(u)z+f''(u)z^2=0, \qquad
z(0)=\sup_{\R}\Big\{\max\big(0,(u_o)_x\big)\Big\} \in \bar \R.
$$
Therefore one recovers (\ref{milo_marat}) in the case $C=1/\kappa$ and $\omega
\equiv Lip(g)$ since we have:
$$
z(t)=\frac{e^{Lip(g).t}}{\frac{1}{z(0)}+\kappa\left(\frac{e^{Lip(g).t}-1}{Lip(g)}\right)}
\leq \frac{e^{Lip(g).t}}{\kappa t}, \qquad t>0.
$$

All in all, the $\lu(\R)$ bound on $\omega$ expresses somehow the fact that
the source has negligible effects outside a compact interval in $\R$ as
pointed out in \cite{daf} p.329. Hence, by the strict hyperbolicity
ensured by (\ref{nonres})--(\ref{nana}), waves exit this
region after some time and then are ruled only by the convective process.

\section*{Acknowledgments}

\noindent
{Both authors acknowledge the hospitality of the Mathematics Department at University of Pavia
where this reseach has been conducted. We also thank Prof. A.~Bressan for constructive comments.}

\bibliographystyle{amsplain}

\end{document}